\documentclass[a4paper,11pt]{amsart}
\usepackage{graphicx}
 \usepackage{mathptmx}
\usepackage{amsmath}
\usepackage{amssymb}
\usepackage{enumitem}
\usepackage{xcolor}

\newmuskip\pFqmuskip

\newcommand*\pFq[6][8]{%
  \begingroup % only local assignments
  \pFqmuskip=#1mu\relax
  \mathcode`=\string"8000
  \begingroup\lccode`\~=`\,
  \lowercase{\endgroup\let~}\pFqcomma
  F^{#2}_{#3}{\left(\genfrac..{0pt}{}{#4}{#5}\bigg|#6\right)}%
  \endgroup
}
\newcommand{\pFqcomma}{\mskip\pFqmuskip}

\newtheorem{theorem}{Theorem}

\begin{document}

\title[A Note on Degenerate Gamma Random Variables]{A Note on Degenerate Gamma Random Variables}

\author{Taekyun  Kim}
\address{School of Science, Xi’an Technological University, Xi’an 710021, China\\
Department of Mathematics, Kwangwoon University, Seoul 139-701, Republic of Korea}
\email{tkkim@kw.ac.kr}

\author{Dae San Kim}
\address{Department of Mathematics, Sogang University, Seoul, Republic of Korea}
\email{dskim@sogang.ac.kr}

\author{Jongkyum Kwon}
\address{Department of Mathematics Education, Gyeongsang National University, Jinju, Republic of Korea}
\email{mathkjk26@gnu.ac.kr}

\author{Hyunseok Lee}
\address{Department of Mathematics, Kwangwoon University, Seoul 139-701, Republic of Korea}
\email{luciasconstant@kw.ac.kr}

\subjclass[2010]{11B73; 11S80; 60G50}
\keywords{degenerate gamma random variable; degenerate gamma fuction; expectation; variance }

\maketitle

\begin{abstract}
In this paper, we introduce the degenerate gamma random variables which are connected with the degenerate gamma functions and the degenerate exponential functions, and deduce the expectation and variance of those random variables.
\end{abstract}

\section{Introduction}
For $0 \neq \lambda\in\mathbb{R}$, the degenerate exponential functions are defined by
\begin{equation}
	e_{\lambda}^{x}(t)=\sum_{n=0}^{\infty}(x)_{n,\lambda}\frac{t^{n}}{n!},\quad e_{\lambda}(t)=e_{\lambda}^{1}(t)=\sum_{n=0}^{\infty}(1)_{n,\lambda}\frac{t^{n}}{n!},\label{1}
\end{equation}
where $(x)_{0,\lambda}=1,\ (x)_{n,\lambda}=x(x-\lambda)\cdots(x-(n-1)\lambda)$, $(n\ge 1)$, (see [5,7,8,9,10]).\\
Note that $\displaystyle \lim_{\lambda\rightarrow 0}e_{\lambda}^{x}(t)=e^{xt}\displaystyle$. \\
As is known, the $\lambda$-binomial coefficients are defined as
\begin{equation}
	\binom{x}{n}_{\lambda}=\frac{(x)_{n,\lambda}}{n!}=\frac{x(x-\lambda)\cdots(x-(n-1)\lambda)}{n!},\quad(n\ge 1),\quad \binom{x}{0}_{\lambda}=1,\quad(\mathrm{see}\ [6,7]).\label{2}
\end{equation}
The degenerate Stirling numbers of the first kind are defined by Kim-Kim as
\begin{equation}
	(x)_{n}=\sum_{l=0}^{n}S_{1,\lambda}(n,l)(x)_{l,\lambda},\quad (n\ge 0),\quad (\mathrm{see}\ [10]),\label{3}
\end{equation}
where $(x)_{0}=1$, $(x)_{n}=x(x-1)\cdots(x-n+1)$, $(n\ge 1)$. \\
As an inversion formula of $\eqref{3}$, the degenerate Stirling numbers of the second kind are given by
\begin{equation}
	(x)_{n,\lambda}=\sum_{l=0}^{n}S_{2,\lambda}(n,l)(x)_{l},\quad(n\ge 0),\quad(\mathrm{see}\ [4]).\label{4}
\end{equation}
Recently, Kim-Kim introduced the degenerate gamma functions which are given by
\begin{equation}
	\Gamma_{\lambda}(s)=\int_{0}^{\infty}e_{\lambda}^{-1}(t)t^{s-1}dt,\quad\big(\lambda\in(0,1)\big),\quad(\mathrm{see}\ [5,9]), \label{5}
\end{equation}
where $s\in\mathbb{C}$, with $0<\mathrm{Re}(s)<\frac{1}{\lambda}$. \\
From \eqref{5}, we note that
\begin{equation}
	\Gamma_{\lambda}(k)=\frac{\Gamma(k)}{(1)_{k+1,\lambda}},\quad\bigg(k\in\mathbb{N},\ \lambda\ne 1,\frac{1}{2},\frac{1}{3},\dots,\frac{1}{k}\bigg), \quad \Gamma_{\lambda}(1)=\frac{1}{1-\lambda}, \quad(\mathrm{see}\ [9]),\label{6}
\end{equation}
where $\Gamma$ is the usual gamma function. \\
Note that $\displaystyle\lim_{\lambda\rightarrow 0+}\Gamma_{\lambda}(s)=\Gamma(s)=\int_{0}^{\infty}e^{-t}t^{s-1}dt, \quad (\textnormal{Re}(s) > 0)\displaystyle$. \\
The random variables are real valued functions defined on sample spaces. We say that $X$ is a continuous random variable if there exists a nonnegative function $f$, defined for all real $x\in(-\infty,\infty)$, having the property that, for any set $B$ of real numbers
\begin{equation}
	P\{X\in B\}=\int_{B}f(x)dx,\quad(\mathrm{see}\ [14]).\label{7}
\end{equation}
The function $f$ is called the probability density function of $X$. If $X$ is a continuous random variable having the probability density function $f$, then the expectation of $X$ is defined as
\begin{equation}
	E[X]=\int_{-\infty}^{\infty}xf(x)dx. \label{8}
\end{equation}
Let $X$ be a continuous random variable with the probability density function $f$. Then, for any real valued function $g$, we have
\begin{equation}
	E[g(X)]=\int_{-\infty}^{\infty}g(x)f(x)dx. \label{9}
\end{equation}
The expected value of the random variable $X$, $E[X]$, is also referred to as the mean or the first moment of $X$. The quantity $E[X^{n}]$, $n\ge 1$, is said to be the $n$-th moment of $X$. That is,
\begin{equation}
	E[X^{n}]=\int_{-\infty}^{\infty}x^{n}f(x)dx,\quad(n\ge 1).\label{10}
\end{equation}
Another quantity of interest is the variance of the random variable $X$ which is defined by
\begin{equation}
	\mathrm{Var}(X)=E\big[(X-E[X])^{2}\big]=E[X^{2}]-\big(E[X]\big)^{2}.\label{11}
\end{equation}
We say that the random variables $X$ and $Y$ are jointly continuous if there exists a function $f(x,y)$, defined for all $x$ and $y$, having the probability that for all sets $A$ and $B$ of real numbers
\begin{equation}
	P\{X\in A,\ Y\in B\}=\int_{B}\int_{A}f(x,y)dxdy.\label{12}
\end{equation}
The function $f(x,y)$ is called the joint probability density function of $X$ and $Y$.
The random variables $X$ and $Y$ are said to be independent if, for all $a,b$, 
\begin{displaymath}
	P\{X\le a,\ Y\le b\}=P\{X\le a\}\cdot P\{Y\le b\},
\end{displaymath}
that is to say, if 
\begin{displaymath}
	E[XY]=E[X]E[Y].
\end{displaymath}
Let $X,Y$ be independent random variables. For any real valued functions $h$ and $g$, we have
\begin{displaymath}
	E[g(X)h(Y)]=E[g(X)]E[h(Y)].
\end{displaymath}
A continuous random variable $X$ is said to have a gamma distribution with parameters $\alpha>0$ and $\beta>0$ if its probability  density function has the form
\begin{displaymath}
	f(x)=\left\{\begin{array}{ccc}
	\frac{1}{\Gamma(\alpha)}\beta e^{-\beta x}(\beta x)^{\alpha-1}, & \textrm{if $x \geq 0,$} \\
	0, & \textrm{otherwise.}
	\end{array}\right.
\end{displaymath}
In this case, we shall say $X$ is a gamma random variable with parameters $\alpha$ and $\beta$, for which we indicate by $X\sim\Gamma(\alpha,\beta)$.
The gamma random variables have long been studied by many researchers (see [1-3, 11-13,16]).
Note that
\begin{align*}
	P\{X\in(-\infty,\infty)\}\ &=\ \int_{0}^{\infty}f(x)dx=\frac{\beta}{\Gamma(\alpha)}\int_{0}^{\infty}(\beta x)^{\alpha-1}e^{-\beta x}dx \\
	&=\ \frac{1}{\Gamma(\alpha)}\int_{0}^{\infty}y^{\alpha-1}e^{-y}dy=1.
\end{align*}
In this paper, we study the degenerate gamma random variables with parameters $\alpha$ and $\beta$ arising from the degenerate gamma function and deduce the expectation and variance of those random variables.

\section{Degenerate Gamma Random Variables}
For $\lambda\in (0,1)$, a random variable $X=X_{\lambda}$ is said to have a degenerate gamma distribution with parameters\, $\alpha$ and $\beta$, \,\,$(\frac{1}{\lambda}>\alpha>0,\,\, \beta > 0)$, \, if its probability density function has the form
\begin{equation}
	f_{\lambda}(x)=\left\{\begin{array}{ccc}
		\frac{1}{\Gamma_{\lambda}(\alpha)}\beta e_{\lambda}^{-1}(\beta x)(\beta x)^{\alpha-1}, & \textrm{if $x>0,$}\\
		0, & \textrm{otherwise.}
	\end{array}\right.\label{13}
\end{equation}
In this case, we shall say that $X$ is a degenerate gamma random variable with parameters $\alpha$ and $\beta$, for which we write as $X\sim \Gamma_{\lambda}(\alpha,\beta)$. \\
Not that
\begin{align*}
	P\{X\in(-\infty,\infty)\}\ &=\ \int_{-\infty}^{\infty}f_{\lambda}(x)dx=\frac{\beta}{\Gamma_{\lambda}(\alpha)}\int_{0}^{\infty}e_{\lambda}^{-1}(\beta x)(\beta x)^{\alpha-1}dx \\
	&=\ \frac{1}{\Gamma_{\lambda}(\alpha)}\int_{0}^{\infty}e_{\lambda}(y)^{-1}y^{\alpha-1}dy=\frac{1}{\Gamma_{\lambda}(\alpha)}\Gamma_{\lambda}(\alpha)=1.
\end{align*}
From \eqref{5}, we note that
\begin{equation}
	\Gamma_{\lambda}(s+1)=\frac{s}{1-\lambda(s+1)}\Gamma_{\lambda}(s).\label{14}
\end{equation}
Let $X\sim\Gamma_{\lambda}(\alpha,\beta)$. For $n\in\mathbb{N}$, the $n$-th moment of $X$ is given by
\begin{align}
	E[X^{n}]\ &=\ \int_{-\infty}^{\infty}x^{n}f_{\lambda}(x)dx\ =\ \frac{\beta}{\Gamma_{\lambda}(\alpha)}\int_{0}^{\infty}e_{\lambda}^{-1}(\beta x)(\beta x)^{\alpha-1}x^{n}dx \label{15} \\
	&=\ \frac{1}{\Gamma_{\lambda}(\alpha)}\frac{1}{\beta^{n}}\int_{0}^{\infty}e_{\lambda}^{-1}(y)y^{n+\alpha-1}dy\ =\ \frac{1}{\beta^{n}}\frac{\Gamma_{\lambda}(n+\alpha)}{\Gamma_{\lambda}(\alpha)}. \nonumber
\end{align}
From \eqref{14}, we can derive the following equation.
\begin{align}
	\Gamma_{\lambda}(n+\alpha)\ &=\ \frac{(n+\alpha-1)\Gamma_{\lambda}(n+\alpha-1)}{1-\lambda(n+\alpha)} \label{16} \\
	&=\ \frac{(n+\alpha-1)(n+\alpha-2)}{\big(1-\lambda(n+\alpha)\big)\big(1-\lambda(n+\alpha-1)\big)}\Gamma_{\lambda}(n+\alpha-2)\\
& =\ \cdots\nonumber \\
& =\ \frac{(n+\alpha-1)(n+\alpha-2)\cdots\alpha\Gamma_{\lambda}(\alpha)}{\big(1-\lambda(n+\alpha)\big)\big(1-\lambda(n+\alpha-1)\big)\cdots\big(1-\lambda(\alpha+1)\big)}.\nonumber
\end{align}
By \eqref{15} and \eqref{16}, we get
\begin{align}
	E[X^{n}]&=\frac{1}{\beta^{n}}\frac{(n+\alpha-1)_{n}}{\big(1-\lambda(\alpha+1)\big)_{n,\lambda}}\label{17}\\
	&=\frac{1}{\beta^{n}}\frac{\binom{n+\alpha-1}{n}}{\binom{1-\lambda(\alpha+1)}{n}_{\lambda}}.\nonumber
\end{align}
Therefore, by \eqref{17}, we obtain the following theorem.
\begin{theorem}
	Let $X\sim\Gamma_{\lambda}(\alpha,\beta)$. For $n\in\mathbb{N}$, we have
	\begin{displaymath}
		E[X^{n}]=\frac{1}{\beta^{n}}\cdot\frac{\binom{n+\alpha-1}{n}}{\binom{1-\lambda(\alpha+1)}{n}_{\lambda}}.
	\end{displaymath}
\end{theorem}
From Theorem 1, we note that
\begin{align*}
	&E[X]\ =\ \frac{\alpha}{\beta\big(1-\lambda(\alpha+1)\big)},\\
	&E[X^2]=\ \frac{1}{\beta^{2}}\frac{\binom{\alpha+1}{2}}{\binom{1-\lambda(\alpha+1)}{2}_{\lambda}}\ =\ \frac{1}{\beta^{2}}\bigg(\frac{\alpha(\alpha+1)}{\big(1-\lambda(\alpha+1)\big)\big(1-\lambda(\alpha+2)\big)}\bigg).
\end{align*}
Thus, the variance of $X$ is given by
\begin{align}
	\mathrm{Var}(X)\ &=\ E[X^{2}]-\big(E[X]\big)^{2} \nonumber \\
	&=\ \frac{1}{\beta^{2}}\bigg(\frac{\alpha(\alpha+1)}{\big(1-\lambda(\alpha+1)\big)\big(1-\lambda(\alpha+2)\big)}\bigg)-\frac{\alpha^{2}}{\beta^{2}}\bigg(\frac{1}{1-\lambda(\alpha+1)}\bigg)^{2} \label{18} \\
	&=\ \frac{\alpha}{\beta^{2}}\bigg(\frac{(\alpha+1)\big(1-\lambda(\alpha+1)\big)-\alpha\big(1-\lambda(\alpha+2)\big)}{\big(1-\lambda(\alpha+1)^{2}\big(1-\lambda(\alpha+2)\big)}\bigg) \nonumber \\
	&=\ \frac{\alpha}{\beta^{2}}\bigg(\frac{1-\lambda}{\big(1-\lambda(\alpha+1)\big)^{2}\big(1-\lambda(\alpha+2)\big)}\bigg). \nonumber
\end{align}
Therefore, by \eqref{18}, we obtain the following theorem
\begin{theorem}
	Let $X\sim \Gamma_{\lambda}(\alpha,\beta)$. Then we have
	\begin{displaymath}
		\mathrm{Var}(X)= \frac{\alpha}{\beta^{2}}\bigg(\frac{1-\lambda}{\big(1-\lambda(\alpha+1)\big)^{2}\big(1-\lambda(\alpha+2)\big)}\bigg).
	\end{displaymath}
\end{theorem}

\section{Further Remark}
Let $X\sim\Gamma_{\lambda}(1,1)$. Then we have
\begin{align}
	E[e_{\lambda}(Xt)]\ &=(1-\lambda)\ \int_{0}^{\infty}e_{\lambda}(xt)e_{\lambda}^{-1}(x)dx \nonumber \\
	&=(1-\lambda)\ \sum_{n=0}^{\infty}(1)_{n,\lambda}\int_{0}^{\infty}x^{n}e_{\lambda}^{-1}(x)dx\frac{t^{n}}{n!}\label{19} \\
	&=\ (1-\lambda) \sum_{n=0}^{\infty}(1)_{n,\lambda}\Gamma_{\lambda}(n+1)\frac{t^{n}}{n!}\nonumber \\
& =\ (1-\lambda)\sum_{n=0}^{\infty}(1)_{n,\lambda}\frac{\Gamma(n+1)}{(1)_{n+2,\lambda}}\frac{t^{n}}{n!}\nonumber \\
	&=\ \sum_{n=0}^{\infty}\frac{1-\lambda}{(1-n\lambda)(1-(n+1)\lambda)}t^{n}.\nonumber
\end{align}
By \eqref{19}, we get
\begin{displaymath}
	(1)_{n,\lambda}\frac{E[X^{n}]}{n!}=\frac{1-\lambda}{(1-n\lambda)(1-(n+1)\lambda)},\quad(n\in\mathbb{N}).
\end{displaymath}
Let $X_{1}\sim\Gamma_{\lambda}(\alpha_{1},\beta_{1})$, $X_{2}\sim\Gamma_{\lambda}(\alpha_{2},\beta_{2})$, $\dots$, $X_{r}\sim\Gamma_{\lambda}(\alpha_{r},\beta_{r})$. If $X_{1},X_{2},\dots,X_{r}$ are identically independent, then we have
\begin{displaymath}
	E[(X_{1}+\cdots+X_{r})^{n}]=\frac{1}{\beta^{n}}\sum_{l_{1}+\cdots+l_{r}=n}\binom{n}{l_{1},\dots,l_{r}}\prod_{i=1}^{r}\binom{l_{i}+\alpha_{i}-1}{l_{i}}\bigg/\binom{1-\lambda(\alpha_{i}+1)}{l_{i}}_{\lambda},
\end{displaymath}
where $n$ is a positive integer. \\
Let $X\sim\Gamma_{\lambda}(\alpha,1)$. Then we have
\begin{align}
	E[e_{\lambda}^{X}(t)]\ &=\ \sum_{n=0}^{\infty}\frac{t^{n}}{n!}E[(X)_{n,\lambda}] \label{20}\\
	&=\ \sum_{n=0}^{\infty}\sum_{l=0}^{n}S_{2,\lambda}(n,l)E[(X)_{l}]\frac{t^{n}}{n!} \nonumber \\
	&=\ \sum_{n=0}^{\infty}\bigg(\sum_{l=0}^{n}\sum_{m=0}^{l}S_{2,\lambda}(n,l)S_{1}(l,m)E[X^{m}]\bigg)\frac{t^{n}}{n!}\nonumber \\
	&=\sum_{n=0}^{\infty}\bigg(\sum_{l=0}^{n}\sum_{m=0}^{l}S_{2,\lambda}(n,l)S_{1}(l,m)\binom{m+\alpha-1}{m}\bigg/\binom{1-(\alpha+1)\lambda}{m}_{\lambda}\bigg)\frac{t^{n}}{n!}. \nonumber
\end{align}
Thus, by \eqref{20}, we get
\begin{equation}
	E[(1+t)^{X}]\ =\ E[e_{\lambda}^{X}\big(\log_{\lambda}(1+t)\big)] \label{21}
\end{equation}
\begin{align*}
	&=\sum_{n=0}^{\infty}\bigg(\sum_{l=0}^{n}\sum_{m=0}^{l}S_{2,\lambda}(m,l)S_{1}(l,m)\binom{m+\alpha-1}{m}\bigg/\binom{1-(\alpha+1)\lambda}{m}_{\lambda}\bigg)\frac{1}{n!}\big(\log_{\lambda}(1+t)\big)^{n} \\
	&=\sum_{n=0}^{\infty}\bigg(\sum_{l=0}^{n}\sum_{m=0}^{l}S_{2,\lambda}(m,l)S_{1}(l,m)\binom{m+\alpha-1}{m}\bigg/\binom{1-(\alpha+1)\lambda}{m}_{\lambda}\bigg)\sum_{k=n}^{\infty}S_{1,\lambda}(k,n)\frac{t^{k}}{k!}\\
	&=\sum_{k=0}^{\infty}\bigg(\sum_{n=0}^{k}\sum_{l=0}^{n}\sum_{m=0}^{l}S_{2,\lambda}(n,l)S_{1}(l,m)S_{1,\lambda}(k,n) \binom{m+\alpha-1}{m}\bigg/\binom{1-(\alpha+1)\lambda}{m}_{\lambda}\bigg)\frac{t^{k}}{k!}.
\end{align*}
On the other hand,
\begin{align}
	E[(1+t)^{X}]\ &=\ \sum_{n=0}^{\infty}E[X^{n}]\frac{1}{n!}\big(\log(1+t)\big)^{n} \label{22} \\
	&=\ \sum_{n=0}^{\infty}E[X^{n}]\sum_{k=n}^{\infty}S_{1}(k,n)\frac{t^{k}}{k!} \nonumber \\
	&=\ \sum_{k=0}^{\infty}\bigg(\sum_{n=0}^{k}S_{1}(k,n)E[X^{n}]\bigg)\frac{t^{k}}{k!}\nonumber \\
	&=\ \sum_{k=0}^{\infty}\bigg(\sum_{n=0}^{k}S_{1}(k,n)\binom{n+\alpha-1}{n}\bigg/\binom{1-(\alpha+1)\lambda}{n}_{\lambda}\bigg)\frac{t^{k}}{k!}\nonumber.
\end{align}
Therefore, by \eqref{21} and \eqref{22}, we obtain the following theorem:
\begin{align*}
	&\sum_{n=0}^{k}S_{1}(k,n)\binom{n+\alpha-1}{n}\bigg/\binom{1-(\alpha+1)\lambda}{n}_{\lambda}\\
	&\ = \sum_{n=0}^{k}\sum_{l=0}^{n}\sum_{m=0}^{l}S_{2,\lambda}(n,l)S_{1}(l,m)S_{1,\lambda}(k,n) \binom{m+\alpha-1}{m}\bigg/\binom{1-(\alpha+1)\lambda}{m}_{\lambda},
\end{align*}
where $X\sim\Gamma_{\lambda}(\alpha,1)$.

\end{document}